\documentclass [11pt,oneside,a4paper,mathscr]{amsart}

\usepackage{amscd,amsmath,amssymb,euscript}
\usepackage[frame,cmtip,curve,arrow,matrix,line,graph]{xy}

\title{An introduction to motivic Hall algebras}
\author{Tom Bridgeland}

\date{}
\newtheorem{thm}{Theorem}[section]

\newtheorem{prop}[thm]{Proposition}
\newtheorem{lemma}[thm]{Lemma}
\newenvironment{pf}{\paragraph{Proof}}{\qed\par\medskip}
\theoremstyle{definition}
\newtheorem{defn}[thm]{Definition}

\newtheorem{thm*}[thm]{Theorem$^*$}
\newtheorem{remark}[thm]{Remark}

\newcommand {\comment}[1]{}

\renewcommand{\leq}{\leqslant}
\renewcommand{\geq}{\geqslant}
\newcommand{\isom}{\cong}
\newcommand{\tensor}{\otimes}
\newcommand{\onto}{\twoheadrightarrow}
\newcommand{\blob}{{\scriptscriptstyle\bullet}}
\newcommand{\into}{\hookrightarrow}

\newcommand{\lRa}[1]{\xrightarrow{\ #1\ }}

\newcommand{\lra}{\longrightarrow}

\newcommand{\dual}{\vee}

\newcommand{\QQ}{\mathbb{Q}}

\newcommand{\can}{\operatorname{can}}

\newcommand{\D}{{D}}

\newcommand{\PP}{\operatorname{\mathbb P}}
\newcommand{\C}{\mathbb C}

\newcommand{\A}{\mathcal A}

\newcommand{\OO}{\mathcal O}

\newcommand{\LL}{\mathbb{L}}

\newcommand{\M}{\mathcal{M}}

\newcommand{\Z}{\mathbb{Z}}
\newcommand{\R}{\mathbf R}

\newcommand{\lExt}{\operatorname{Ext}}
\newcommand{\Coh}{\operatorname{Coh}}

\newcommand{\id}{\operatorname{id}}
\newcommand{\Ext}{\operatorname{Ext}}
\newcommand{\Hom}{\operatorname{Hom}}
\newcommand{\Map}{\operatorname{Map}}
\newcommand{\Spec}{\operatorname{Spec}}

\newcommand{\GL}{\operatorname{GL}}

\newcommand{\St}{\operatorname{St}}
\newcommand{\virt}{\operatorname{vir}}

\newcommand{\var}{\operatorname{Var}}
\newcommand{\red}{{\operatorname{red}}}

\newcommand{\RHom}{\operatorname{\mathbf{R} Hom}}
\newcommand{\RlHom}{\operatorname{\mathbf{R}\mathcal{H} om}}

\newcommand{\sch}{\operatorname{Sch}}
\renewcommand{\sp}{\operatorname{Sp}}

\newcommand{\Isom}{\operatorname{Isom}}

\newcommand{\an}{\operatorname{an}}

\newcommand{\ksch}[1]{K({\operatorname{Sch/{#1}}})}
\newcommand{\kvar}[1]{K({\operatorname{Var/{#1}}})}

\newcommand{\ksp}[1]{K({\operatorname{Sp/{#1}}})}
\newcommand{\kst}[1]{K({\operatorname{St/{#1}}})}
\newcommand{\RH}{\operatorname{H}}

\newcommand{\RHreg}{\operatorname{H_{{reg}}}}
\newcommand{\RHsc}{\operatorname{H_{sc}}}

\newcommand{\affsch}[1]{\operatorname{Aff/{#1}}}

\begin{document}

\begin{abstract}
We give an introduction to Joyce's construction of the motivic
Hall algebra of coherent sheaves on a variety $M$. When $M$ is a
Calabi-Yau threefold we define a semi-classical integration map
 from a Poisson subalgebra of this Hall algebra to the ring of functions on a
symplectic torus.  This material will be used in \cite{forth} to
prove some basic properties of Donaldson-Thomas curve-counting
invariants on Calabi-Yau threefolds.
\end{abstract}
\maketitle

\section{Introduction}

This paper is a gentle introduction to part of Joyce's theory of
motivic Hall algebras \cite{Jo0,Jo1,Jo2,JS}. It started life as
the first half of the author's paper \cite{forth} in which this
theory is used to prove some properties of Donaldson-Thomas
curve-counting invariants on Calabi-Yau threefolds. Eventually it
became clear that there
 were enough points at which our presentation differs from Joyce's  to justify a separate paper.
  Nonetheless, most of the basic ideas  can  be found in
Joyce's work.
%\medskip

The application of Hall algebras to the study  of invariants of
moduli  spaces  originated with Reineke's computation of the Betti
numbers of the  spaces of stable quiver representations
\cite{rei}. His technique was to   translate
categorical statements %(for example, existence and uniqueness of
%Harder-Narasimhan filtrations)
 into identities in a suitable Hall
algebra, and to then apply a ring homomorphism into a completed
skew-polynomial ring, thus obtaining identities involving
generating functions for the invariants of interest.
 %\medskip

The relevant category  in Reineke's paper is the category of
representations of a finite quiver without relations. Such
categories can be defined over any field, and Reineke worked with
a Hall algebra based on counting points over $\mathbb{F}_q$. In
\cite{Jo2} Joyce used Grothendieck rings  of Artin stacks to
construct a motivic version of the Hall algebra defined in
arbitrary characteristic. This can be applied, for example, to
categories of coherent sheaves on complex varieties.

%\medskip

The interesting part of the theory is the construction of a
homomorphism from the  Hall algebra to a skew polynomial ring,
often viewed as a ring of functions on a quantum torus.  Such maps
go under the general name of integration maps, since they involve
integrating an element of the Hall algebra over the moduli stack.
In Reineke's case the existence of such a map relied  on the fact
that the relevant categories  of representations were of
homological dimension one. Remarkably it seems that integration
maps also exist when the underlying abelian category
 is Calabi-Yau of dimension 3.

 %\medskip
 In the CY$_3$ case Joyce's integration map is a homomorphism of Lie algebras defined on a Lie subalgebra
 of the Hall algebra.
Kontsevich and Soibelman \cite{KS} suggested that incorporating motivic vanishing cycles should enable one to construct an algebra morphism from the full Hall algebra.
 Unfortunately the details of the constructions in their paper are currently  rather sketchy  and rely  on
some unproved conjectures.
%\medskip
Joyce and Song \cite{JS} went on to use some of the ideas from
\cite{KS} to prove an important property of Behrend functions
(stated here as Theorem \ref{see}) and so incorporate such
functions into Joyce's Lie algebra integration map.

The main result of this paper (Theorem \ref{se}) is the existence
of an integration map in the CY$_3$ case that is a homomorphism of
Poisson algebras. It can be viewed as the semi-classical limit  of
the ring homomorphism envisaged by Kontsevich and Soibelman. It
relies on the same property of the Behrend function proved by
Joyce and Song.  Combined with a difficult no-poles result of
Joyce it can be used to prove non-trivial results on
Donaldson-Thomas invariants \cite{forth}.
\subsection*{Acknowledgements} Thanks most of all to Dominic Joyce who  patiently
explained many things about motivic Hall algebras. Thanks also to
Arend Bayer, Andrew Kresch and  Max Lieblich for useful
conversations. Finally, I'm very grateful to Maxim Kontsevich and
Yan Soibelman for sharing a preliminary version of their paper
\cite{KS}.

\renewcommand{\H}{\mathcal{H}}

% ******************************************************************************************************
% ******************************************************************************************************
% ******************************************************************************************************

\section{Grothendieck rings of varieties and schemes}

Here we review some basic definitions concerning Grothendieck
rings of varieties. Some good references are \cite{bitt,loo}. For
us a complex variety is a reduced, separated scheme of finite type
over $\C$. We denote by
\begin{equation}\label{blbl}\var/\C\subset \sch/\C\subset\sp/\C\end{equation}
the categories of varieties over $\C$, of schemes of finite type over $\C$,
 and of algebraic spaces of finite type over $\C$ respectively.
%We also consider categories
%\[\sch_\infty/\C\subset\sp_\infty/\C\]
%of schemes and algebraic spaces locally of finite type over $\C$.

\subsection{Grothendieck ring of varieties}
\label{fir} Recall first the definition of the Grothendieck ring
of varieties.

\begin{defn}\label{scc}Let $\kvar{\C}$ denote the free abelian group on isomorphism classes
 of complex varieties, modulo relations
\begin{equation}
\label{scissor}
[X]=[Z]+[U]\end{equation} for $Z\subset X$ a closed subvariety with
complementary open subvariety $U$.
\end{defn}

The relations \eqref{scissor} are called the scissor relations,
since they involve cutting a variety up into pieces. We can equip
$\kvar{\C}$ with the structure of a commutative ring by setting
\[[X]\cdot [Y]=[X\times Y].\]
The class of a point $1=[\Spec(\C)]$ is then a unit. We write
\[\LL=[\mathbb{A}^1] \in \kvar{\C}\]
for the class of the affine line.
%\medskip
By a stratification of a variety $X$ we mean a  collection of
disjoint locally-closed subsets $X_i\subset X$ which together
cover $X$.

\begin{lemma}
\label{silk} If a variety $X$ is stratified by subvarieties $X_i$
then only finitely many of the $X_i$ are non-empty and
\[[X]=\sum_i [X_i]\in \kvar{\C}.\]
\end{lemma}

\begin{pf}The result is clear for varieties of dimension 0, so let
us use induction on the dimension $d$ of $X$, and assume the
result known for varieties of dimension $<d$.

  First consider the
case when $X$ is irreducible. Then one of the $X_i=U$ contains the
generic point and is therefore open. The complement $Z=X\setminus
U$ is of smaller dimension and is stratified by the other
subvarieties $X_i$. Since \[[X]=[Z]+[U]\] the result follows by
induction.

 In the case when $X$ is reducible we can take an
irreducible subvariety and remove the intersections with the other
irreducible components. This gives  an irreducible open subset
$U\subset X$ with complement a closed subvariety $Z$ having fewer
irreducible components than $X$. By induction on this number one
can therefore conclude that
\[[Z]=\sum_i[Z\cap X_i], \quad [U]=\sum_i [U\cap X_i], \]
with finitely many non-empty terms appearing in each sum. Since
\[[X_i]=[Z\cap X_i]+[U\cap X_i],\]
the result then follows from the scissor relations.
\end{pf}

 There is a ring homomorphism $\chi\colon \kvar{\C}\to
\Z$ defined by sending the class of a variety $X$ to its topological Euler characteristic
\[\chi(X)=\sum_{i=0}^{2d} (-1)^i \dim H^i(X_{\an},\C), \]
where $X_{\an}$ denotes $X$ equipped with the analytic topology, and $H^i$ denotes singular cohomology. %\smallskip

\begin{lemma}
\label{grimer}
Suppose $x\in\kvar{\C}$ satisfies \[\LL^m\cdot (\LL^n-1) \cdot x =0\] for some $m,n\geq 1$. Then $\chi(x)=0$.
\end{lemma}

\begin{pf}
There is a ring homomorphism  $\chi_t\colon \kvar{\C}\to \Z[t]$
that sends the class of a smooth complete variety $X$ to the
Poincar{\'e} polynomial
\[\chi_t(X)=\sum_{i=0}^{2d} t^i \cdot \dim H^i(X_{an},\C).\]
It specializes at $t=-1$ to the Euler characteristic.
Now
\[\chi_t(\LL)=\chi_t(\PP^1)-\chi_t(1)=t^2.\]
Since $\Z[t]$ is an integral domain  one therefore has $\chi_t(x)=0$.
Setting $t=-1$ gives the result.
\end{pf}

% ******************************************************************************************************
% ******************************************************************************************************
% ******************************************************************************************************

\subsection{Zariski fibrations}
 There is a useful identity in $\kvar{\C}$ relating to
 fibrations.

\begin{defn}A morphism of schemes $f\colon X\to Y$ will be called  a Zariski fibration if there is an open cover $Y=\bigcup_{i\in I} U_i$ and  diagrams
\begin{equation*}\xymatrix@C=1em{
f^{-1}(U_i)\ar[rr]^{g_i} \ar[dr]_{f}&& U_i\times
F_i\ar[dl]^{\pi_1}\\  &U_i }
\end{equation*}with each $g_i$ an isomorphism.

 \end{defn}

Of course if $Y$ is connected then all the fibres $F_i$ are
 isomorphic, but this will often not be the case.
We say that two Zariski
fibrations
 \[f_1 \colon X_1 \to Y\text{  and  }f_2 \colon X_2 \to Y\]
 have the same fibres, if for any point $y\in Y(\C)$ the fibres of $f_1$ and $f_2$ over $y$  are isomorphic.

\begin{lemma}
\label{fib} Suppose
 $f_1\colon X_1\to Y$ and $f_2\colon X_2\to Y$ are
Zariski  fibrations of varieties with the same fibres. Then \[[X_1]=[X_2]\in
\kvar{\C}.\]
\end{lemma}

\begin{pf}
We can stratify $Y$ by a finite collection of connected, locally-closed subvarieties $Y_i\subset Y$ such that
  $f_1$ and $f_2$ are trivial fibrations over each $Y_i$.
   Then
  \[[X_1]=\sum_i [f_1^{-1}(Y_i)]=\sum_i [F_i]\cdot [Y_i]=\sum_i [f_2^{-1}(Y_i)]=[X_2],\]
  where $F_i$ is the common fibre of $f_1$ and $f_2$ over $Y_i$.
  \end{pf}

The following application of Lemma \ref{fib}  will be important
later.

\begin{lemma}
There is an identity
\[[\GL_d]=\LL^{ \frac{1}{2} d (d-1)}\cdot  \prod_{k=1}^{d} (\LL^{k} -1)\in
\kvar{\C}.\]
\end{lemma}

\begin{pf}
Let $B\subset \GL_d$
be the stabilizer of a nonzero vector
 $x\in \C^d$. The assignment $g\mapsto g(x)$ defines a morphism \[\pi\colon \GL_d\to\C^d\setminus\{0\}\]
 which is easily seen to be a Zariski fibration with fibre  $B$.
 But there is also an isomorphism \[B\isom \GL_{d-1}\times\C^{d-1}.\]
   Thus by Lemma \ref{fib} \[[\GL_d]=(\LL^{d}
-1)\cdot \LL^{d-1} \cdot [\GL_{d-1}],\] and the result follows by
induction.
\end{pf}

% ******************************************************************************************************
% ******************************************************************************************************
% ******************************************************************************************************

\subsection{Geometric bijections}
We will base our treatment of Grothendieck groups on the following
class of maps.

\begin{defn}
\label{gb} A morphism $f\colon X \to Y$ in the category $\sch/\C$
is a geometric bijection if it induces a bijection \[f(\C) \colon
X(\C) \to Y(\C)\] between the sets of $\C$-valued points.
\end{defn}

 Using Lemma \ref{bij} below, it is not difficult to prove
that the condition of Definition \ref{gb} is equivalent to $f$
being a bijection, or a universal bijection, but for various
reasons we prefer to introduce new terminology.

By a stratification of a scheme $X$ we mean a collection of
disjoint locally-closed subschemes $X_i \subset X$ which together
cover $X$. If the scheme $X$ is of finite type over $\C$ then the
argument of Lemma \ref{silk} shows that  only finitely many of the
$X_i$ can be non-empty.
%It is a universal bijection if it remains a bijection after pulling back via any morphism of schemes.
%such that the morphism
%\[\coprod_\alpha i_\alpha\colon \coprod_\alpha Y_\alpha \to Y \]
%is a bijection. This is equivalent to the statement that the  %Recall that the image of any immersion is  locally-closed.

\begin{lemma}
\label{bij} A morphism $f\colon X\to Y$  in the category $\sch/\C$
is a geometric bijection precisely if there
 are stratifications
\[X_i\subset X, \quad Y_i\subset Y, \]
 such that $f$ induces
isomorphisms $f_i\colon X_i \to Y_i$.
\end{lemma}

\begin{pf}A more precise statement of the condition is that there
should be isomorphisms $f_i\colon X_i \to Y_i$ and commuting
diagrams
\[\begin{CD} X_i &@>f_i>> &Y_i \\
@Vj_iVV && @VVk_iV \\
X &@>f>> &Y \\\end{CD}%\medskip
\]
where the morphisms $j_i$ and $k_i$ are the embeddings of the
given locally-closed subschemes. This condition is clearly
sufficient since every $\C$-valued point of $X$ or $Y$  factors
through a unique one of the given subschemes.

For the converse we may as well assume that $X$ and $Y$ are
reduced, since a stratification of $X_{\red}$ also gives a
stratification of $X$, and similarly for $Y$.  We claim that there
is an open subscheme $Y_1\subset Y$ such that $f$ induces an
isomorphism
\[f\colon X_1\to Y_1, \]
where $X_1=f^{-1}(Y_1)$.
 This will be enough since we can then replace $X$ and
$Y$ by the  complements of $X_1$ and $Y_1$ and
repeat.

 To prove the claim we can pass to an  open subset of $Y$ and hence assume that
$Y$ is an irreducible variety. By generic flatness, we can also
assume that $f$ is flat, and hence open. Now we can replace $X$ by
an irreducible open subvariety, and so $f$ becomes a map of
irreducible varieties.  The claim then holds by \cite[Prop.
3.17]{mum}.
 \end{pf}

Using Lemma \ref{bij} we can  give an alternative definition of
the Grothendieck ring in terms of bijections. This is sometimes
useful, particularly when considering Grothendieck rings of
schemes and stacks of possibly infinite type.

\begin{lemma}\label{sccc}The  group $\kvar{\C}$ is the free abelian group on isomorphism classes
 of the category $\var/\C$, modulo relations
\begin{itemize}
\item[(a)] $[X_1\amalg X_2]=[X_1]+[X_2]$ for every pair of
varieties $X_1$ and $X_2$,%\smallskip

 \item[(b)] $[X]=[Y]$ for
every geometric bijection $f\colon X\to Y$.
\end{itemize}
\end{lemma}

\begin{pf}
Given a geometric bijection of varieties $f\colon X \to Y$, we can
take stratifications of $X$ and $Y$ as in Lemma \ref{bij}. We can
always assume that the subschemes $X_i$ and $Y_i$ are reduced and
hence subvarieties.  Then,  by Lemma \ref{silk},
\[[X]=\sum_i [X_i]=\sum_i [Y_i] = [Y] \in \kvar{\C}.\]
Thus relation (b) is a consequence of the scissor relations, and
clearly relation (a) is a special case of them. Conversely, given
a decomposition  as in Definition \ref{scc}, the obvious morphism
$Z\amalg U\to X$ is  a geometric bijection, so the scissor
relations of Definition \ref{scc} are a consequence of the
relations (a) and (b) in the statement of the Lemma.
\end{pf}

% *************************************************************************************************************************
% *************************************************************************************************************************
% *************************************************************************************************************************

\subsection{Grothendieck rings of schemes and algebraic spaces}
\label{inf}

Before considering stacks in the next section,
 it is worth briefly  considering the case  of schemes and algebraic spaces, always of finite type over $\C$.

\begin{defn}
\label{mantle}
Let $\ksch{\C}$ be the free abelian group on isomorphism classes
 of the category $\sch/\C$, modulo relations
\begin{itemize}
\item[(a)] $[X_1\amalg X_2]=[X_1]+[X_2]$ for every pair of
schemes $X_1$ and $X_2$,%\smallskip

 \item[(b)] $[X]=[Y]$ for
every geometric bijection $f\colon X\to Y$.
\end{itemize}
\end{defn}

The product in $\sch/\C$ gives the group $\ksch{\C}$ the structure
of a commutative ring. One could alternatively define $\ksch{\C}$
via scissor relations as in Definition \ref{scc}; the argument of
Lemma \ref{sccc} shows that this would give the same ring.

In the case of algebraic spaces we  define geometric bijections in
the category $\sp/\C$ exactly as in Definition \ref{gb}. We can
also define the notion of a
stratification of an algebraic space in the obvious way.% The
%analogue of Lemma \ref{bij}

\begin{lemma}
\label{bijspace} A morphism $f\colon X\to Y$  in the category
$\sp/\C$ is a geometric bijection precisely if there
 are stratifications
\[X_i\subset X, \quad Y_i\subset Y, \]
 such that $f$ induces
isomorphisms $f_i\colon X_i \to Y_i$.
\end{lemma}

\begin{pf}
This follows the same lines as that of Lemma \ref{bij}. The extra
argument needed is the following. Suppose given a morphism
$f\colon X\to Y$ in $\sp/\C$ with $X$ and $Y$ reduced. We must
show that there is an open subspace $Y_1\subset Y$ such that $f$
induces an isomorphism \[f\colon X_1\to Y_1,\]
 where $X_1=f^{-1}(Y_1)$.
By \cite[Prop. II.6.6]{Knutson} we can pass to an open subset and
so assume that $Y$ is a scheme, and even an irreducible variety.
By generic flatness we can also assume that $f$ is flat. Then,
using the same result from \cite{Knutson} again, there is an open
subset $X_0\subset X$ representable by an irreducible variety.
Since the induced map $f\colon X_0\to Y$ is flat, its image is an
open subvariety $Y_0\subset Y$. We can then apply \cite[Prop.
3.17]{mum} as in the proof of Lemma \ref{bij}.
\end{pf}

The Grothendieck group $\ksp{\C}$ is defined as in Definition
\ref{mantle}, replacing the category $\sch/\C$ by $\sp/\C$. The
following result shows that from the point of view of Grothendieck
rings, providing we stick to objects of finite type, the
distinction between varieties, schemes and algebraic spaces
disappears.

 \begin{lemma}
 \label{proofabove}
 The embeddings of categories \eqref{blbl} induce isomorphisms of rings \[\kvar{\C}\isom\ksch{\C}\isom
 \ksp{\C}.\]
 \end{lemma}

\begin{pf}
The basic point is that if $Y\in \sp/\C$  there is a geometric
bijection
\[f\colon X\to Y\]
with $X$ a variety. Indeed, by \cite[Prop. II.6.6]{Knutson}  there
is an open subspace $U\subset Y$
  that is representable by an affine scheme.
 Taking the complement and repeating, we can stratify $Y$ by affine schemes $Y_i \subset Y$.
  The inclusion map from the disjoint union of these strata then defines a geometric bijection \[X=\coprod Y_i \to Y\]
  with $X$ an affine scheme of finite type over $\C$. We can
 assume that $X$ is reduced and hence a variety
   since the inclusion of its reduced subscheme is another geometric bijection.

 Now consider the
homomorphism
\[I\colon \kvar{\C} \to \ksp{\C}\]
induced by the inclusion of varieties in algebraic spaces. By the
above it is surjective, so it will be enough to construct a left
inverse
\[P \colon \ksp{\C} \to \kvar{\C}\]
Given an object $Y\in \sp/\C$, take a geometric bijection $f\colon
X\to Y$ with $X$ a variety, and set \[P([Y])=[X].\] This is
well-defined, because if $W$ is another variety with a geometric
bijection $g\colon W\to Y$, then we can form the fibre square
\[\begin{CD} Z &@>p>> &W \\
@VqVV && @VVgV \\
X &@>f>> &Y \end{CD}
\]
and taking a variety $T$ with a geometric bijection $g\colon T\to
Z$, the composite morphisms $p\circ g$ and $q\circ g$ are
geometric bijections, so \[[X]=[W]\in \kvar{\C}.\] It is easy to
check that $P$ preserves the relations and hence defines the
required inverse.
\end{pf}

% ******************************************************************************************************
% ******************************************************************************************************
% ******************************************************************************************************

\section{Grothendieck rings of stacks}

In this section we consider Grothendieck rings of algebraic
stacks. The author learnt most of the material in this section
from the papers of Joyce \cite{Jo0} and To{\"e}n \cite{Toen}. We
refer to \cite{G} for a readable introduction to stacks, and to
\cite{LMB} for a more detailed treatment. All stacks will be
Artin\comment{Quasi-separated condition?} stacks and will be
assumed to be locally of finite type over $\C$.
 We denote by
$\St/\C$ the 2-category of algebraic stacks of finite type over
$\C$.

Given a scheme $S$ and a stack $X$ we denote by $X(S)$ the
groupoid of $S$-valued points of $X$. We will be fairly sloppy
with 2-terminology: by a commutative (resp. Cartesian) diagram of
stacks we mean one that is 2-commutative (resp. 2-Cartesian). By
an isomorphism of stacks we mean what a category-theorist would
call an equivalence.

% ******************************************************************************************************
% ******************************************************************************************************
% ******************************************************************************************************

\subsection{Geometric bijections and Zariski fibrations}

In the case of stacks the appropriate analogue of Definition
\ref{gb} is as follows.

\begin{defn}
\label{blblbl}
 A morphism $f\colon X\to Y$ in the category $\St/\C$ will be called a geometric bijection
  if  it is representable and the induced functor on groupoids of $\C$-valued points
\[f(\C)\colon X(\C)\to Y(\C)\]
is an equivalence of categories.
\end{defn}

In fact, it is easy enough\comment{Add more explanation} to show
that the assumption that $f$ be representable in Definition
\ref{blblbl} follows from the other condition. This comes down to
showing that if $f\colon X \to Y$ is a group scheme in $\sch/\C$
such that for each point $y\in Y(\C)$ the fibre $X_y$ is the
trivial group, then $f$ is an isomorphism.

By a stratification of a stack $X$
 we mean a collection of locally-closed substacks $X_i\subset X$
  that are disjoint and together cover $X$.
 % By this last condition we mean that given an atlas $Y\to X$
%  the induced collection of subschemes $Y_i\subset Y$ cover $Y$.
If $X$ is of finite type over $\C$ it follows by pulling back to
an atlas that only finitely many of the $X_i$ can be non-empty.

\begin{lemma}
\label{stratastack}
 A morphism $f\colon X\to Y$  in
$\St/\C$ is a geometric bijection precisely if there are
stratifications
\[X_i\subset X, \quad Y_i\subset Y, \]
 such that $f$ induces
isomorphisms $f_i\colon X_i \to Y_i$.
\end{lemma}

\begin{pf}
As in the proof of Lemma \ref{bij} it is easy to see that the
condition is sufficient. For the converse, assume that $f$ is a
geometric bijection. Replacing $X$ and $Y$ by their reduced
substacks we can assume that $X$ and $Y$ are reduced. It will be
enough to show that there is a non-empty open substack $Y_1
\subset Y$ such that $f$ induces an isomorphism
\[f\colon X_1 \to Y_1,\] where $X_1=f^{-1}(Y_1)$. We can then take complements and repeat.

Pulling $f$ back to an atlas $\pi\colon T\to Y$ we obtain
a diagram \[\begin{CD} S &@>g>> &T \\
@VVV && @VV\pi V \\
X &@>f>> &Y \\\end{CD}
\]with $S$ an algebraic space and $g$ a
geometric bijection.

Let $T_1\subset T$ be the largest open subset over which $g$ is an
isomorphism.  By Lemma \ref{bijspace}
 this subset is non-empty.  That it descends to $Y$ follows from the following statement.
 Suppose given a Cartesian
 diagram in $\sp/\C$
\begin{equation}
\label{diaa}\begin{CD} S' &@>g'>> &T'\\
@VqVV && @VVpV \\
S &@>g>> &T \\\end{CD}
\end{equation}
 with $p$ faithfully flat. Then  $g'$ is an isomorphism over an open subset $U'\subset T'$
 precisely if $g$ is an isomorphism over the open subset $U=p(U')\subset
 T$.

  Replacing $T$ by $U$ and $T'$ by $U'$ we can reduce further to the statement that given a
 diagram \eqref{diaa} with $p$ faithfully flat and $g'$ an
 isomorphism then $g$ is also an isomorphism.
 This is the well-known statement that isomorphisms are stable in
 the faithfully flat topology \cite[Prop. 3.5]{Knutson}.
\end{pf}

\begin{defn}A morphism of stacks $f\colon X\to Y$ is a
Zariski fibration if its pullback to any scheme is a Zariski fibration of schemes.
 %We say that two Zariski fibrations over a stack have the same fibres if this is true when pulled back to any atlas.
\end{defn}

  In particular a Zariski fibration of stacks is  representable.
 Note however that there need not be a cover of the stack $Y$ by
  open substacks over which $f$ is trivial. For this reason, in the Grothendieck group of stacks,
   the fibration
  identity
  of Lemma \ref{fib} is not a consequence of the other relations, and we
 will impose it by hand.

%\medskip
% ******************************************************************************************************
% ******************************************************************************************************
% ******************************************************************************************************

\subsection{Grothendieck ring of stacks}

To make the comparison result in  the next section true, it will
be necessary to restrict slightly the class of  stacks we
allow in our Grothendieck ring. %Recall that an algebraic group is affine
 %if it is isomorphic to a closed subgroup of $\GL_d$.

\begin{defn}
A stack $X$ locally of finite type over $\C$  has affine
stabilizers if for every $\C$-valued point $x \in X(\C)$ the group
$\Isom_\C(x,x)$ is affine.
\end{defn}

 The importance of this
notion lies in the following corollary of a result of Kresch.

\begin{prop}
\label{kresch}
 A stack $X\in \St/\C$ has affine  stabilizers precisely if there is a variety $Y$ with an action of $G=\GL_d$ and
 a geometric bijection
 \[f\colon Y/G\to X.\]
\end{prop}

\begin{pf}
The condition is clearly sufficient. For the converse, suppose
first that the automorphism groups of all geometric points of $X$
are affine. Kresch then shows \cite[Prop. 3.5.2, Prop. 3.5.9]{K}
that
 the associated reduced stack $X_{\red}$ has a
stratification by locally-closed substacks of the form
\[X_i=Y_i/G_i,\] where by \cite[Lemma 3.5.1]{K}
we can take each group $G_i$ to be of the form $G_i=\GL_{d_i}$.
  The obvious
map $f\colon \amalg X_i \to X$ is then a geometric bijection, and
the result follows from the isomorphism
\[(X_1/G_1) \amalg (X_2/G_2) \isom [(X_1\times G_2) \sqcup (X_2\times
G_1)]/(G_1\times G_2),\] which shows that the disjoint union of
quotient stacks is another quotient stack.

To finish the proof we must check that if the automorphism groups
of  all $\C$-valued points of $X$ are affine, then the same is
true for all geometric points. To prove this,\footnote{The author
is grateful to Andrew Kresch for explaining this argument.}
suppose $f\colon G \to S$ is a group scheme, with $G$ and $S$ of
finite type over $\C$, and suppose that there is a geometric point
$\Spec(K)\to S$ such that the corresponding geometric fibre of $f$
is non-affine. We must prove that there is a $\C$-valued point of
$S$ such that the corresponding fibre of $f$ is non-affine.

 Since being affine
is invariant under field extension [EGA IV.2.7.1] we can assume
that $K$ is the algebraic closure\comment{Why need algebraic
closure?} of the residue field $k(s)$ of a point $s\in S$.
Chevalley's Theorem \cite{Chev} implies that an algebraic group
$G$ over a
 field $K$ of characteristic zero is
non-affine precisely if there is an epimorphism $G_0 \onto A$
where $G_0\subset G$ is the connected component of the identity,
and $A$ is a positive-dimensional abelian variety. For a modern
proof of this see \cite{Con}.

Applying this to our situation, the epimorphism $G_0\onto A$ is
defined over a finite extension of $k(s)$, and
hence\comment{Explain this} over a finite type scheme with a
morphism $T\to S$ dominating the closure of the point $s$.
Restricting to a $\C$-valued point and applying Chevalley's
theorem again completes the proof.
\end{pf}

We can now give the following definition\comment{Sort out $X_1$,
$X_2$ versus $X$ and $Y$}.

\begin{defn}
\label{mainy}Let $\kst{\C}$ be the free abelian group spanned by
isomorphism classes of
  stacks of finite type over $\C$ with affine
 stabilizers, modulo   relations
\begin{itemize}
\item[(a)] $[X_1\amalg X_2]=[X_1]+[X_2]$ for every pair of stacks
$X_1$ and $X_2$,

%\smallskip
 \item[(b)] $[X]=[Y]$ for every geometric bijection
$f\colon X\to Y$,
%\smallskip

\item[(c)] $[X_1]=[X_2]$ for every pair of Zariski fibrations
$f_i\colon X_i \to Y$ with the same fibres.
\end{itemize}
\end{defn}

Fibre product of stacks over $\C$ gives $\kst{\C}$ the structure
of a commutative ring. There is an obvious homomorphism of
commutative rings
\begin{equation}
\label{obvious}\kvar{\C}\to\kst{\C}\end{equation} obtained by
considering a variety as a representable stack.

% ******************************************************************************************************
% ******************************************************************************************************
% ******************************************************************************************************

\subsection{Comparison Lemma}

 In this section, following To{\"e}n \cite{Toen}, we show  that if one localises $\kvar{\C}$ at the classes of all special algebraic
groups, the map \eqref{obvious} becomes an isomorphism.

\begin{defn}An  algebraic group is
 special if any map of schemes \[f\colon X\to Y\] which is a principal $G$-bundle  in the
  {\'e}tale topology is a Zariski fibration.
 \end{defn}

Examples of special groups include the general linear groups
$\GL_d$. All special groups are affine and connected \cite{Chow}.

\begin{lemma}
Localizing the ring $\kvar{\C}$ with respect to any of the following three sets of elements gives the same result:
\begin{itemize}
\item[(a)]  the classes $[G]$ for $G$ a special algebraic
group,%\smallskip

\item[(b)]  the classes $[\GL_d]$ for $d\geq 1$,%\smallskip
 \item[(c)]  the elements $\LL$
and $\LL^i-1$ for $i\geq 1$.
\end{itemize}
\end{lemma}

\begin{pf}
 If $G\subset \GL_d$ is a closed subgroup,
then the quotient map \[\pi\colon \GL_d \to \GL_d/G\]
 is a
principal $G$-bundle. It is locally trivial in the smooth topology
(equivalently, in the {\'e}tale topology \cite[Example 2.51]{vis})
because it is smooth, and becomes trivial when pulled back to
itself. Thus if $G$ is special we can conclude that
\[[\GL_d]=[G]\cdot [\GL_d/G].\] Since $\GL_d$ is itself special
this proves the equivalence of (a) and (b). The equivalence of (b)
and (c) follows from Lemma \ref{fib}.
\end{pf}

If an algebraic group $G$ acts on a variety $X$ then the map to
the quotient stack
\[\pi\colon X \to X/G\]
is a principal $G$-bundle when pulled back to any scheme. If $G$
is special it follows that $\pi$ is a Zariski fibration with fibre
$G$ and hence \begin{equation} \label{matrix}[X]=[G]\cdot
[X/G]\in\kst{\C}.\end{equation} In particular taking $X$ to be a
point we see that $[G]$ is
invertible in $\kst{\C}$. %\medskip

The following result is due to To{\"e}n \cite[Theorem 3.10]{Toen} (who also considered higher stacks).
A slightly weaker version  was proved independently by Joyce
 \cite[Theorem 4.10]{Jo0}.

\begin{lemma}
\label{toen} The homomorphism \eqref{obvious} induces an isomorphism of commutative rings
\[Q\colon \kvar{\C}[[\GL(d)]^{-1}:d\geq 1] \lra \kst{\C}.\]
\end{lemma}

\begin{pf}
The expression \eqref{matrix} shows that the map $Q$ is
well-defined and satisfies
\[Q([X]/[G])=[X/G]\in\kst{\C}.\]
We shall construct an inverse to $Q$.
 Suppose $Z$ is a stack with affine
stabilizers. By Proposition \ref{kresch} there is a geometric
bijection $g\colon X/G\to Z$ with $G$ a special algebraic group
acting on a variety $X$. Set
\[R(Z)=[X]/[G] \in \kvar{\C}[[\GL(d)]^{-1}:d\geq 1]. \]
To check that this is well-defined, suppose given another
bijection $h\colon Y/H\to Z$. Consider the diagram of Cartesian
squares
\[\begin{CD}
R@>p>> &Q  @>j>> &Y \\
@VqVV & @VVV  &@VVV \\
P @>>> &W @>>> &Y/H\\
@VkVV & @VVV  &@VVh V \\
X @>\pi>> &X/G @>g>> &Z \\ \end{CD}\] Then $Q$ is an algebraic
space and $j$ is a geometric bijection, so $[Q]=[Y]\in \ksp{\C}$.
The map $\pi$, and hence also $p$, is a Zariski fibration with
fibre $G$; pulling it back to a bijection $Q'\to Q$ with $Q'$ a
scheme, it follows that $[R]=[G]\cdot [Q]\in\ksp{\C}$. Hence we
obtain $[R]=[G]\cdot [Y]$ and, by symmetry, $[R]=[H]\cdot [X]$.
Thus, using Lemma \ref{proofabove}
\[[G]\cdot [Y]=[H]\cdot [X] \in \kvar{\C}.\]

To show that $R$ descends to the level of the Grothendieck group
we must check that the relations in Definition \ref{mainy} are
mapped to zero. This is very easy for (a) and (b). For (c) suppose
$f_i\colon X_i\to Y$ are Zariski fibrations with the same fibres.
Take a geometric bijection $W/G\to Y$ and form the diagrams
\[
\begin{CD}
S_i @>k_i>> &T_i @>h_i>> &X_i\\
@Vp_iVV & @VVV  &@VVf_iV \\
W @>\pi>> &W/G @>g>> &Y  \end{CD}\]
 Then there are induced actions of $G$ on the varieties $S_i$ such that $T_i\isom S_i/G$. Since $h_i$ is a geometric bijection,
 $R(X_i)=[S_i]/[G]$. On the other hand, by pullback, the morphisms $p_i\colon S_i\to W$ are Zariski fibrations
 of schemes with the same fibres, and hence by the argument of Lemma \ref{fib}, $[S_1]=[S_2]\in \ksch{\C}$.
\end{pf}

%It follows from Lemma \ref{toen} that any motivic invariant for
%varieties   for which $\Upsilon(\GL(d))$ is
%invertible for all $d$ induces a ring homomorphism
%\begin{equation*} \Upsilon\colon \kst{\C}\to
%\Lambda.\end{equation*} In particular, since  $\chi_t(\LL)=t^2$ it
%follows that $\chi_t$
% induces a ring homomorphism
% \[\chi_t\colon \kst{\C}\to \Z(t).\]
%

% ******************************************************************************************************
% ******************************************************************************************************
% ******************************************************************************************************

\subsection{Relative Grothendieck groups}

Let $S$ be a fixed algebraic stack, locally of finite type over
$\C$. We shall always assume that $S$ has affine stabilizers. There
is a 2-category of algebraic stacks over $S$.
 Let $\operatorname{St/S}$ denote the full subcategory  consisting of objects
\begin{equation}
\label{obvious2}f \colon X \to S\end{equation} for which $X$ is of
finite type over $\C$. Such an object will be said to have affine
stabilizers if the stack $X$ has. Repeating Definition \ref{mainy}
in this relative context gives the following.
%\begin{equation} \label{sym} [X\lRa{f}S],\end{equation}
%where $X$ is a stack of finite type and $f$ is a
%morphism. Two such symbols
%\[[X_1\lRa{f_1}S]\quad\text{ and }\quad[X_2\lRa{f_2}S]\] will
%be called equivalent if there is a commutative diagram
%with $g$ an isomorphism.

\begin{defn}
\label{rel} Let $\kst{S}$ be the free abelian group spanned by
isomorphism classes of objects \eqref{obvious2} of
$\operatorname{St/S}$,
 with affine stabilizers,  modulo relations
  \begin{itemize}
\item[(a)] for every pair of objects $X_1$ and $X_2$ a
relation\[[X_1\amalg X_2\lRa{f_1\sqcup f_2} S]=[X_1\lRa{f_1}
S]+[X_2\lRa{f_2} S],\]

\item[(b)] for every commutative diagram
\[ \xymatrix@C=.8em{
X_1\ar[dr]_{f_1}\ar[rr]^{g} && X_2\ar[dl]^{f_2}\\
&S }
\]with $g$ a geometric bijection, a  relation
\[[X_1\lRa{f_1} S]=[X_2\lRa{f_2} S],\]

\item[(c)]  for every pair of Zariski fibrations
%\[ \xymatrix@C=.8em{
%X_1\ar[dr]_{h_1}&& X_2\ar[dl]^{h_2}\\
%&Y }
%\]
\[h_1\colon X_1 \to Y, \quad h_2\colon X_2\to Y\]
  with the same fibres, and
every morphism $g\colon Y\to S$, a relation
\[[X_1\lRa{g\circ h_1} S]=[X_2\lRa{g\circ h_2} S].\]
\end{itemize}
\end{defn}

The group $\kst{S}$ has the structure of a $\kst{\C}$-module,
defined by setting
\[[X]\cdot [Y\lRa{f} S]=[X\times Y\lRa{f\circ\pi_2} S]\]
and extending linearly.

\begin{remark}
Suppose that $\Lambda$ is a $\QQ$-algebra and
\[\Upsilon \colon \kst{\C}\lra \Lambda\]
is a ring homomorphism. Then for each stack $S$ with affine
stabilizers Joyce defines \cite[Section 4.3]{Jo2} a
$\Lambda$-module $\operatorname{SF}(S,\Upsilon,\Lambda)$ whose
elements he calls stack functions. It is easy to see that there is
an isomorphism of $\Lambda$-modules
\[\operatorname{SF}(S,\Upsilon,\Lambda)\isom \kst{S}\tensor_{\kst{\C}}
\Lambda.\] We leave the  proof to the reader.
\end{remark}

% ******************************************************************************************************
% ******************************************************************************************************
% ******************************************************************************************************

\subsection{Functoriality}
\label{funct}
 The following statements are are all easy consequences of
the basic properties of fibre products of stacks, and we leave the
proofs to the reader. We assume that all stacks appearing  have affine
stabilizers.

 %\medskip

\begin{itemize}
%\smallskip

\item[(a)]A morphism of stacks $a\colon S \to T$ induces a map of
$\kst{\C}$-modules
\[ a_* \colon \kst{S} \lra \kst{T}\]
sending $[X\lRa{f} S]$ to $[X\lRa{a\circ f} T]$.

\item[(b)] A morphism of stacks $a\colon S \to T$ of finite type
 induces a map of
$\kst{\C}$-modules
\[a^*\colon \kst{T}\lra \kst{S}\]
sending $[ Y\lRa{g} T]$
to %the pullback
%\[f\times_T S \colon  X\times_T S \lra S.\] %
 $[X\lRa{f} S]$ where $f$ is the map appearing in the Cartesian
square
\[\begin{CD}
X& @>>> &Y\\
 @VfVV  &&@VVg V \\
S &@>a>> &T\end{CD}\]

\item[(c)]The above assignments are functorial, in that
\[(b\circ a)_* = b_* \circ a_*, \quad (b\circ a)^* = a^*\circ b^*,\]
whenever $a$ and $b$ are composable morphisms of stacks with the
required properties.

 \item[(d)]Given a Cartesian square of maps
\[\begin{CD} U &@>c>>& V\\
 @VdVV && @VVbV \\
S& @>a>> &T\end{CD}\] one has the base-change property \[b^*\circ
a_* = c_*\circ d^*\colon \kst{S}\lra \kst{V}.\]

\item[(e)] For every pair of stacks $(S_1,S_2)$ there is a
K{\"u}nneth map
\[K\colon \kst{S_1} \tensor \kst{S_2} \to \kst{S_1 \times S_2}\]
given by
\[[X_1\lRa{f_1} S_1] \tensor [X_2\lRa{f_2} S_2] \mapsto
[X_1\times X_2 \lRa{f_1\times f_2} S_1\times S_2].\]
It is a morphism of $\kst{\C}$-modules. %\smallskip
\end{itemize}

We can  view the functor $\kst{--}$ as defining a primitive
cohomology theory for stacks.

% ******************************************************************************************************
% ******************************************************************************************************
% ******************************************************************************************************

\section{Motivic Hall algebra}
\label{hall} Let $M$ be a smooth projective variety and
\[\A=\Coh(M)\] its category of coherent sheaves.
In this section,
following Joyce \cite{Jo2}, we introduce the motivic Hall algebra
of the category $\A$. Much of what we do here
would apply with minor modifications to other abelian categories,  %(for example categories of finite-dimensional modules over an associative algebra),
but we make no attempt at maximal generality.

We use the following abuse of notation throughout: if $f\colon T\to S$ is a morphism
of schemes, and $E$ is a sheaf on $S\times M$, we use the shorthand $f^*(E)$ for the pullback to $T\times M$, rather than the more correct $(f\times 1_M)^*(E)$.

% ******************************************************************************************************
% ******************************************************************************************************
% ******************************************************************************************************

\subsection{Stacks of flags}

Let $\M^{(n)}$ denote the moduli stack of $n$-flags of coherent
sheaves on $M$. The objects over a scheme $S$ are chains of
monomorphisms of coherent sheaves on $S\times M$ of the form
\begin{equation}
\label{star} 0=E_0\into E_1 \into \cdots \into E_n=E\end{equation}
such that each factor $F_i=E_i/E_{i-1}$ is $S$-flat. It follows
that each sheaf $E_i$ is also $S$-flat. If
\[0=E'_0\into E'_1 \into \cdots \into E'_n=E\]
is another such object over a scheme $T$, then a morphism in
$\M^{(n)}$ lying over a morphism of schemes $f\colon T\to S$ is a
collection of isomorphisms of sheaves
\[\theta_i \colon f^*(E_i) \to E'_i\]
such that each diagram
\[\begin{CD} f^*(E_i) &@>>> &f^*(E_{i+1})\\
 @V\theta_iVV && @VV\theta_{i+1}V \\
E'_i& @>>> &E'_{i+1}\end{CD}\] commutes. Here we have taken the
usual step of choosing, for each morphism of sheaves, a pullback
of every coherent sheaf on its target. The stack property for
$\M^{(n)}$ follows easily from the corresponding property of the
stack $\M =\M^{(1)}$.
%\medskip

There are morphisms of stacks
\[a_i\colon \M^{(n)}\to \M, \quad 1\leq i\leq n, \]
sending a flag \eqref{star} to its $i$th factor $F_i=E_i/E_{i-1}$.
To define these it is first necessary to choose a cokernel for
each monomorphism $E_{i-1} \to E_i$. There is another morphism
\[b\colon \M^{(n)}\to \M\]
sending a flag \eqref{star} to the sheaf $E_n=E$. Note that the
functors defining all these morphisms of stacks have the iso-fibration property of Lemma \ref{headache}.
 Using this it is immediate that for $n>1$ there is a Cartesian
 square
\begin{equation}\label{first}\begin{CD}
\M^{(n)} & @>t>> &\M^{(2)} \\
 @VsVV & &@VVa_1V \\
\M^{(n-1)} &@>b>> &\M\end{CD}\end{equation} where $s$ and $t$ send
a flag \eqref{star} to the flags \[E_1\into \cdots\into
E_{n-1}\text{ and }E_{n-1}\into E_n\]
 respectively.
%\medskip

There is a kind of duality around here, which is the basic reason
for the associativity of the Hall algebra. Instead of considering
flags of the form \eqref{star} we could instead consider flags
\begin{equation} \label{dagger} E=E^0\onto E^{1}\onto \cdots \onto
E^{n-1}\onto E^n=0.\end{equation} Setting $E^i=E/E_i$ shows that
this gives an isomorphic stack. This dual approach leads to
Cartesian diagrams
\begin{equation}\label{second}\begin{CD}
\M^{(n)} & @>v>> &\M^{(2)} \\
 @VuVV & &@VVa_2V \\
\M^{(n-1)} &@>b>> &\M\end{CD}\end{equation} where $u$ and $v$ send
a flag \eqref{dagger} to the flags
\[E^1\onto \cdots \onto E^{n-1}\text{ and }E^{0}\onto E^1\]
 respectively.

%If\[\begin{CD}
%Z & @>f>> &W \\
%@VgVV && @VhVV \\
%X  &@>j>> &Y
%\end{CD}
%\]
%is a Cartesian square then
% \[\begin{CD}
%Z & @>f>> &W \\
%@V(g,k\circ f)VV && @V(h,k)VV \\
%X \times T  &@>(j,1)>> &Y\times T
%\end{CD}
%\]is also.

The stack $\M^{(2)}$ can be thought of as the stack of short exact
sequences in $\A$. There is a diagram
\begin{equation}\label{blib}\begin{CD}
\M^{(2)} &@>b>> \M\\
 @V(a_1,a_2)VV \\
\M\times\M\end{CD}\end{equation}
 where, as above, the morphisms
$a_1,a_2$ and $b$ send a short exact sequence
\[0\lra A_1\lra B\lra A_2\lra 0\]
 to the sheaves $A_1$, $A_2$ and $B$ respectively. %The morphism $b$ is representable but not usually of finite type.
%The fibre over a point of $\M$ corresponding to a sheaf $B$ is
%Grothendieck's Quot scheme
%\[\operatorname{Quot}_M(B)\] parameterizing all quotients of $B$ in the category
%$\A$.

\begin{lemma}
The stacks $\M^{(n)}$ are algebraic. \end{lemma}

\begin{pf}
Suppose $f\colon S \to \M$ is a morphism of stacks corresponding
to a flat family of sheaves $B$ on $S\times M$. Forming the
Cartesian square
\[\begin{CD}
Z & @>>> &\M^{(2)} \\
 @VVV & &@VVbV \\
S &@>f>> &\M\end{CD}\]
 it is easy to see that $Z$ is
represented by the relative Quot scheme parameterising quotients
of $B$ over $S$. Thus $b$ is representable, and pulling back an
atlas for $\M$ gives an atlas for $\M^{(2)}$. Since fibre products
of algebraic stacks are algebraic the result then follows by
induction and  the existence of the squares \eqref{first}.
\end{pf}

The morphism $(a_1,a_2)$ is  not representable. The fibre  over a
point of $\M \times \M$ corresponding to a pair of sheaves
$(A_1,A_2)$ is the quotient stack
\[[\Ext^1(A_2,A_1)/\Hom(A_2,A_1)],\]
with the action of the vector space $\Hom_\A(A_2,A_1)$ being the
trivial one. This statement follows from  Proposition \ref{tom}
below.

\begin{lemma}
The morphism $(a_1,a_2)$ is of finite type.
\end{lemma}

\begin{pf}
Fix a projective embedding of $M$. For  each integer $m$ and each
polynomial $P$  there is a finite type open substack
$\M_m(P)\subset \M$ parameterising $m$-regular sheaves\comment{Give reference} with
Hilbert polynomial $P$. Define an open substack
\[Y=(a_1,a_2)^{-1}(\M_m(P_1)\times \M_m(P_2))\subset \M^{(2)}.\]
It will be enough to show that $Y$ is of finite type. But the
morphism $b$ restricts to a map \[b\colon  Y \to \M_m(P_1+P_2)\]
since an extension of two $m$-regular sheaves is also $m$-regular.
This map is of finite type, because once one fixes the Hilbert
polynomials involved, the relative Quot scheme is of finite type.
\end{pf}

% ******************************************************************************************************
% ******************************************************************************************************
% ******************************************************************************************************

\subsection{The Hall algebra}
%\label{conv}

  Let us set \[\RH(\A)=\kst{\M}.\] Applying the results of Section \ref{funct} to the diagram
    \eqref{blib} gives a morphism of $\kst{\C}$-modules
\[m=b_*\circ (a_1,a_2)^* \colon \RH(\A)\tensor \RH(\A) \lra \RH(\A)\]
which we call the convolution product.\footnote{Here and in the
rest of this section we  suppress an application of the
K{\"u}nneth map
\[\kst{\M}\tensor \kst{\M} \to \kst{\M\times\M}\]
from the notation.}
%\medskip
Explicitly this is given by the rule
\[ [X_1\lRa{f_1}\M] * [X_2\lRa{f_2} \M] = [Z\lRa{b\circ h}\M], \]
where  $h$ is defined by the following Cartesian
square
\[\begin{CD}
Z & @>h>> &\M^{(2)} &@>b>> \M\\
 @VVV  &&@VV(a_1,a_2)V \\
X_1\times X_2 &@>f_1\times f_2>> &\M\times\M\end{CD}\]
%\smallskip

 The following
result is due to Joyce \cite[Theorem 5.2]{Jo2}, although the basic
idea is of course the same as for previous incarnations of the
Hall algebra.

\begin{thm}
The  product $m$ gives $\RH(\A)$ the structure of an associative
unital algebra over $\kst{\C}$. The unit  element is
\[1=[\M_{0}\subset \M],\] where $\M_{0}\isom\Spec(\C)$ is the
substack of zero objects in $\A$.
\end{thm}

\begin{pf}
Consider the composition
\[\begin{CD}
\RH(\A)\tensor \RH(\A)\tensor \RH(\A)@>m\tensor \id>>
\RH(\A)\tensor \RH(\A) @>m>> \RH(\A)\end{CD}.\] It is induced by
the diagram
\[\begin{CD}
 &&&& \M^{(2)} &@>b>> &\M \\
&&&& @V(a_1,a_2)VV \\
\M^{(2)}\times \M &@> (b,\id) >> &\M\times \M \\
@V (a_1,a_2,\id) VV \\
\M\times\M\times \M \end{CD} \] There is a bigger commutative
diagram obtained by filling in the top left square with
\[\begin{CD}
\M^{(3)} & @>t>> &\M^{(2)} \\
@V(s,a_2\circ t)VV && @VV(a_1,a_2)V \\
\M^{(2)} \times \M  &@>(b,\id)>> &\M\times\M
\end{CD}
\]
where $s$ sends a flag $E_1\into E_2\into E_3$ to the flag
$E_1\into E_2$, and $t$ sends it to $E_2\into E_3$. This square is
Cartesian because of the square $(\ref{first})$ and Lemma
\ref{nonsense}, so by the base-change property of Section
\ref{funct}
 \[m\circ (m\tensor \id) = b_* \circ (a_1,a_2,a_3)^*\]
 is induced by the diagram
 \[\begin{CD} \M^{(3)} @>b>> \M \\
 @V(a_1,a_2, a_3)VV  \\
\M^3  \end{CD}\]

A similar argument using the  square $(\ref{second})$ shows that
the other composition is induced by the same diagram. The
multiplication is therefore associative. We leave the reader to
check the unit property.
\end{pf}

\begin{lemma}
The $n$-fold product
\[m_n\colon \RH(\A)^{\tensor n} \to \RH(\A)\]
is induced by the diagram
 \[\begin{CD} \M^{(n)} @>b>> \M \\
 @V(a_1,\cdots, a_n)VV  \\
\M^n \end{CD}\]
in the sense that
\[m_n=b_*\circ (a_1,\cdots,a_n)^* \colon \RH(\A)^{\tensor n}\lra \RH(\A).\]
\end{lemma}

\begin{pf}
This follows by induction and a similar argument to the one above, but using the Cartesian diagram
\[\begin{CD} \M^{(n)} & @>t>> &\M^{(2)} \\
 @V(s,a_2\circ t)VV & &@VV(a_1,a_2)V \\
\M^{(n-1)}\times\M &@>(b,\id)>> &\M\times\M\end{CD}\] given by
\eqref{first} and Lemma \ref{nonsense}.
\end{pf}

% ******************************************************************************************************
% ******************************************************************************************************
% ******************************************************************************************************

\subsection{Grading}

Let $K(M)=K(\A)$ denote the Grothendieck group of the category
$\A$. Given two coherent sheaves $E$ and $F$ we can define
\[\chi(E,F)=\sum_i (-1)^i \dim_{\C}\,\Ext^i(E,F).\]
This defines a bilinear form $\chi(-,-)$ on $K(M)$ called the
Euler form. Serre duality implies that the left and right radicals
$^{\perp}K(M)$ and $K(M)^{\perp}$ are equal. The numerical
Grothendieck group is the quotient
\[N(M)=K(M)/K(M)^{\perp}.\]
Let $\Gamma\subset N(M)$ denote the monoid of effective classes,
which is to say classes of the form $[E]$ with $E$ a sheaf.

\begin{lemma}
\label{stan} Suppose $S$ is a connected scheme and $F$ is an
$S$-flat coherent sheaf on $S\times M$. For each point
$s\in S(\C)$ let
 \[F_s=F|_{\{s\} \times M}\]
 be the corresponding sheaf on $M$. Then the class $[F_s]\in N(M)$ is independent of the point $s$.
\end{lemma}

\begin{pf}
For any locally-free sheaf $E$, the integer
\[\chi(E,F_s)=\chi(E^{\dual}\tensor F_s)\]
is locally constant on $S$. Since $M$ is smooth, the Grothendieck
group $K(M)$ is spanned by the classes of locally-free sheaves,
 so  the class $[F_s]\in N(M)$ is also locally constant.
 \end{pf}

It follows from Lemma \ref{stan} that the stack $\M$ splits as a
disjoint union of open and closed substacks
\[\M=\bigsqcup_{\alpha\in \Gamma} \M_\alpha,\]
where $\M_\alpha\subset \M$ is the stack of objects of class
$\alpha\in \Gamma$.  The inclusion $\M_\alpha\subset \M$ induces
an embeddding
\[\kst{\M_\alpha}\subset \kst{\M}.\]
There is thus a direct sum decomposition
\begin{equation}
\label{grading} \RH(\A)=\bigoplus_{\alpha\in \Gamma}
\RH(\A)_\alpha,
\end{equation}
and $\RH(\A)$ with the convolution product becomes a
$\Gamma$-graded algebra.

\subsection{Sheaves supported in dimension $\leq d$}

For any integer $d$ there is a full abelian subcategory \[\A_{\leq
d}=\Coh_{\leq d}(M)\subset\A=\Coh(M)\] closed under extensions,
consisting of coherent sheaves on $M$ whose support has dimension
$\leq d$. There is a subgroup  \[N_{\leq d}(M) \subset N(M)\]
spanned by classes of objects of $\Coh_{\leq d}(M),$ and a
corresponding positive cone\[\Gamma_{\leq d}=N_{\leq d}(M) \cap
\Gamma.\]
%\medskip
One can define a Hall algebra $\RH(\A_{\leq d})$ by replacing the
moduli stack $\M$ in the above discussion with the substack
$\M_{\leq d}$ of objects of $\A_{\leq d}$.

In fact
 it is easy to see that a coherent sheaf $E$ lies in $\Coh_{\leq
d}(M)$ precisely if its class $[E]\in N(M)$ lies in the subgroup
$N_{\leq d}(M)$. Thus
\[\M_{\leq d} = \bigsqcup_{\alpha\in \Gamma_{\leq d}} \M_\alpha,\]
and  there is an identification
\[\RH(\A_{\leq d}) = \bigoplus_{\alpha\in \Gamma_{\leq d}}
\RH_\alpha(\A).\]
We will make heavy use of the $d=1$ version of this construction in \cite{forth}.

%\begin{lemma}
%There is a ring homomorphism \[R_{\leq d}\colon \RH(\A)\to
%\RH(\A_{\leq d})\] obtained by taking fibre product with the
%inclusion $\M_{\leq d}\subset \M$.
%\end{lemma}
%
%\begin{pf} %lie in $\A_{\leq k}$.
%%This means that there is a Cartesian diagram
%%\[\begin{CD}
%%\M^{(2)}_{\leq k}& @>>> &\M^{(2)}\\
%% @VbVV & &@VVbV \\
%%\M_{\leq k} &@>>> &\M\end{CD}\]
%%The result follows easily from this.
%Left to the reader. The basic point is that given a short exact
%sequence
%\[0\lra A_1\lra B\lra A_2\lra 0\]
%in $\A$, the object $B$  lies in $\A_{\leq d}$ if and only  if the
%same is true of the objects $A_1$ and $A_2$.
%\end{pf}

% ******************************************************************************************************
% ******************************************************************************************************
% ******************************************************************************************************

\section{Integration map}

In this section we construct a homomorphism of Poisson algebras
from a semi-classical limit of the Hall algebra to an algebra of
functions on a symplectic torus. It can be viewed as the
semi-classical limit of the ring homomorphism envisaged by
Kontsevich and Soibelman \cite{KS}. We assume throughout that $M$
is a smooth projective Calabi-Yau threefold over $\C$. We include
in this the condition that \[H^1(M,\OO_M)=0.\] There are two
versions of the story depending on a choice of sign
$\sigma\in\{\pm 1\}$ which we fix throughout. The choice
$\sigma=+1$ will lead to naive Euler characteristic invariants,
while $\sigma=-1$ leads to Donaldson-Thomas invariants.

%\subsection{Behrend function}

\subsection{Regular elements}
Consider  the  maps of commutative rings \[\kvar{\C} \to
\kvar{\C}[\LL^{-1}] \to \kst{\C},\] and recall that $\RH(\A)$ is
an algebra over  $\kst{\C}$. Define a $\kvar{\C}[\LL^{-1}]$-module
\[\RHreg(\A) \subset \RH(\A)\] to be the span of classes of maps
$[X \lRa{f} \M]$ with $X$ a variety. We call an element of
$\RH(\A)$ regular if it lies in this submodule. The following
result will be proved in Section \ref{firstproof} below.

\begin{thm}
\label{fi} The submodule  of regular elements is closed under the
convolution product:
\[\RHreg(\A) * \RHreg(\A) \subset \RHreg(\A),\]
and is therefore a $\kvar{\C}[\LL^{-1}]$-algebra. Moreover the
 quotient
\[\RHsc(\A)=\RHreg(\A)/(\LL-1) \RHreg(\A)\]
is a commutative $\kvar{\C}$-algebra.
\end{thm}

 We call
the algebra $\RHsc(\A)$ the semi-classical Hall algebra. Since
\[[\C^*]=\LL-1\] is invertible in $\kst{\C}$, there is a
Poisson bracket on $\RH(\A)$ given by the formula
\[\{f,g\}=\frac{f*g-g*f}{\LL-1} .\]
This bracket preserves the subalgebra $\RHreg(\A)$ because the
multiplication in $\RHreg(\A)$ is commutative modulo the ideal
$(\LL-1)$. The induced bracket on $\RHreg(\A)$ then descends to
give a Poisson bracket on the commutative algebra $\RHsc(\A)$.

 %and the canonical map
%\[\T_{[t]}[N(X)]\tensor_{\Lambdao} \Lambda \lra \T_{(t)}[N(X)]\]
%is an isomorphism. Again $\T_{[t]}[N(X)]$ is closed under the
%product $*$,

% ******************************************************************************************************
% ******************************************************************************************************
% ******************************************************************************************************

\subsection{The integration map}
Define a ring\[\Z_\sigma[\Gamma]=\bigoplus_{\alpha\in \Gamma}
\Z\cdot x^{\alpha}
\] by taking the free abelian group spanned by symbols
$x^\alpha$ for $\alpha\in\Gamma$ and setting
\[x^\alpha * x^\beta = \sigma^{\chi(\alpha,\beta)} \cdot x^{\alpha +\beta}.\]
Since the Euler form is skew-symmetric this ring is commutative.
Equip it with a Poisson structure  by defining
\[\{x^\alpha,x^\beta\}= \sigma^{\chi(\alpha,\beta)}\cdot
\chi(\alpha,\beta) \cdot x^{\alpha+\beta}.\] Take a locally
constructible function
\[\lambda\colon\M\to\Z.\]
For definitions of constructible functions on stacks see
\cite{Jo-1}. Associated to every sheaf $E\in \A$ is an integral
weight
\[\lambda(E)\in\Z.\]
 If $X$ is a variety with a map $f \colon X \to \M$, there is an induced constructible function
\[f^*(\lambda)\colon X \to \Z,\] and  a
weighted Euler characteristic
\[\chi(X,f^*(\lambda))=\sum_{n\in \Z} n \cdot \chi((\lambda\circ f)^{-1}(n)).\]
  We prove the following result  in Section \ref{secondproof} below.

\begin{thm}
\label{se}Given a locally constructible function $\lambda\colon\M\to\Z$, there is a unique group homomorphism
\[I \colon \RHsc(\A) \lra \Z_\sigma[\Gamma],\]such that if $X$ is a variety and
$f\colon X\to \M$ factors via $\M_\alpha\subset \M$ then
\[I\big([X\lRa{f}\M]\big)=\chi(X,f^*(\lambda)) \cdot x^\alpha.\]
Moreover, $I$ is a homomorphism of commutative algebras if for all
$E,F\in\A$ \begin{equation} \label{firstt}\lambda(E\oplus
F)={\sigma}^{\chi(E,F)} \cdot \lambda(E)\cdot
\lambda(F),\end{equation} and a homomorphism of Poisson algebras
if, in addition, the expression
\begin{equation}
\label{secondd}m(E,F)=\chi(\PP\Ext^1(F,E),
\lambda(G_\theta)-\lambda(G_0))\end{equation} is
symmetric in $E$ and $F$.
\end{thm}

This statement may need a little explanation. We have written $G_\theta$ for
the sheaf
\[0\lra E\lra G_\theta \lra F \lra 0\]
corresponding to a class $\theta\in \Ext^1(F,E)$. In particular $G_0=E\oplus F$.
 For nonzero $\theta$ the isomorphism
class of the object $G_\theta$ depends only on the class
\[[\theta]\in \PP \Ext^1(F,E),\] and it is easy to see that the map $[\theta]\mapsto \lambda(G_\theta)$
 is a constructible function\comment{Prove this} on
the projective space $\PP\Ext^1(F,E)$.

% ********************************************************************************************
% ********************************************************************************************
% ********************************************************************************************

\subsection{Behrend function}

Recall \cite{Be} that Behrend
associates to any scheme $S$ of finite type over $\C$ a
 constructible function
\[\nu_S\colon S\to\Z\]
with the property that if $S$ is a proper\comment{proper =
complete?} moduli scheme with a symmetric obstruction theory then
the associated Donaldson-Thomas virtual count coincides with the
weighted Euler characteristic:
\[\#_{\virt}(S) := \int_{[S]^{\virt}} 1=\chi(S,\nu_S).\]
 These functions satisfy the relation
\begin{equation}
\label{mum} f^*(\nu_S)=(-1)^d \nu_T\end{equation} whenever
$f\colon T\to S$ is a smooth morphism of relative dimension $d$.
It follows easily from this that every stack $S$, locally of
finite type over $\C$, has an associated locally constructible
function
\[\nu_S\colon S\to\Z\] defined uniquely by the condition that \eqref{mum}
 also holds for smooth morphisms of stacks.

 It is clear that the constant function \[\mathbf{1}\colon \M \to \Z\]
 satisfies
the conditions of Theorem \ref{se} with the sign $\sigma=+1$.
Regarding the Behrend function, Joyce and Song \cite[Theorem
5.9]{JS} proved the following wonderful result. For their proof
(which uses gauge-theoretic methods) it is essential that our base
variety $M$ is proper and that we are working over $\C$.

\begin{thm}
\label{see}  The Behrend function
\[\nu_\M\colon \M \to \Z\] for the moduli stack $\M$  satisfies
the conditions of Theorem \ref{se} with the sign $\sigma=-1$.
\end{thm}

 Thus we have (at least) two integration maps: taking $\lambda=\mathbf{1}$ leads to
invariants defined by unweighted Euler characteristics, whereas taking $\lambda=\nu_\M$ leads to Donaldson-Thomas invariants.

%\item[(c)] One can easily eliminate the use of the virtual
%Poincar{\'e} polynomial in this section and throughout the paper,
%although the result is at least to the author less transparent. To
%do this replace the inclusion of rings $\C[t,t^{-1}]\subset \C(t)$
%by the homomorphism
%\[\kvar{\C}[\LL^{-1}] \to \kst{\C}.\]
% Similarly, rather than specialising at the ideal $t-1$
%specialise at $\LL-1$. Then Theorem \ref{fi} still holds, and one
%obtains a Poisson algebra over $\kvar{\C}$. Tensoring with the
%Euler characteristic recovers the algebra $\RH_{-1}(\A)$.

% ******************************************************************************************************
% ******************************************************************************************************
% ******************************************************************************************************

\section{Fibres of the convolution map}
 In this section we give a careful analysis of the fibres of the  morphism of
stacks $(a_1,a_2)$ appearing in the definition of the convolution
product. This will be the main tool in the proofs of Theorems
\ref{fi} and \ref{se} we give in the next section. We again adopt
the abuse of notation for pullbacks of families of sheaves
explained in the preamble to Section 4.

\subsection{Universal extensions}
Suppose $S$ is an affine scheme and $E_1$ and $E_2$ are coherent
sheaves on $S\times M$, flat over $S$. Let $\affsch{S}$ denote the
category of affine schemes over $S$.
 Define a functor
\[\Phi^k_S(E_1,E_2) \colon \affsch{S} \to \operatorname{Ab}\]
by sending an object $f\colon T\to S$  to the
abelian group
\[\Phi^k_S(E_1,E_2)(f)=\Ext^k_{T\times M} (f^*(E_1),f^*(E_2)).\]
The image of  a morphism\[ \xymatrix@C=.8em{
U\ar[dr]_{g}\ar[rr]^{h} && T\ar[dl]^{f}\\
&S }
\]
 in $\affsch{S}$ is defined  using the canonical map
\[h^*\colon \Ext^k_{T\times M} (f^*(E_1),f^*(E_2))\to
\Ext^k_{U\times M}(h^*f^*(E_1),h^*f^*(E_2))\] together with the
canonical isomorphisms
\[\can \colon g^*(E_i) \isom h^* f^* (E_i). \]
To check that this does indeed define a functor one needs to apply
the uniqueness properties of pullback in the usual way (see for
example \cite[Section 3.2.1]{vis}).
%\medskip

Consider the object
\[\RlHom_{\OO_S} (E_1,E_2)=\R \pi_{S,*}\RlHom_{\OO_{S\times M}}(E_1,E_2) \in \D\Coh(S).\]
%In particular, pulling back to a point, we obtain
%\[H^i(V_{12}(f) \Ltensor _{\OO_{T}}
%\OO_t)\isom \Ext^i_X(E_1|_{\{t\}\times X}, E_2|_{\{t\}\times
%X}).\]
For each $k\geq 0$ we set
\[\lExt^k_{\OO_S} (E_1,E_2):=H^k(\RlHom_{\OO_S} (E_1,E_2))\in \Coh(S).\]
 We shall say that $E_1$ and $E_2$
have constant $\Ext$ groups if these sheaves are all locally-free.

\begin{lemma}
\label{stanage} Suppose $E_1$ and $E_2$ are $S$-flat coherent
sheaves on $S\times M$ with constant $\Ext$ groups. Then the
functor \[\Phi^k_S(E_1,E_2) \colon \affsch{S} \to
\operatorname{Ab}\] defined above  is represented by the
vector bundle
 $V^k(S)$
over $S$ corresponding to the locally-free sheaf $\lExt^k_{\OO_S}
(E_1,E_2)$.
\end{lemma}

\begin{pf}
If $f\colon T \to S$ is a morphism of schemes then by flat
base-change\comment{Reference}, \[\RlHom_{\OO_T}
(f^*(E_1),f^*(E_2))\isom \mathbf{L} f^* \circ \RlHom_{\OO_S}
(E_1,E_2).\] If $E_1$ and $E_2$ have constant $\Ext$ groups it
follows that
\[\lExt^k_{\OO_T} (f^*(E_1),f^*(E_2))\isom f^*\lExt^k_{\OO_S} (E_1,E_2).\]
 There is also an identity
\[\RHom_{T\times M}(f^*(E_1),f^*(E_2))\isom \R\Gamma
\circ \RHom_{\OO_T}(f^*(E_1),f^*(E_2)).\] If $T$ is affine it follows
that
\[\Phi^k_S(E_1,E_2)(f)\isom \Gamma(T,f^* (\lExt^k_{\OO_S}
(E_1,E_2)))\isom\Map_S(T,V^k(S)).\] These isomorphisms commute with
pullback\comment{tricky}
 and hence define an isomorphism of functors.
\end{pf}

% ******************************************************************************************************
% ******************************************************************************************************
% ******************************************************************************************************

\subsection{The convolution morphism}

The  main result of this section is as follows.

\begin{prop}
\label{tom} Let $X_1$ and $X_2$ be varieties with morphisms
\[f_1\colon X_1 \to \M,\qquad f_2\colon X_2\to\M,\] and let
$E_i\in \Coh(X_i\times M)$ be the corresponding families of
sheaves on $M$. Then we can stratify $X_1\times X_2$ by
locally-closed affine subvarieties
\[S\subset X_1\times X_2\] with the following property. For each
 point $s\in S(\C)$ the space
\[V^k(s)=\Ext^k_M(E_2|_{\{s\}\times M}, E_1|_{\{s\}\times
M})\] has a fixed dimension $d_k(S)$, and if we form the Cartesian
squares
\[\begin{CD}
Z_S &@>>> &Z & @>h>> &\M^{(2)} \\
@VuVV  && @VtVV  &&@VV(a_1,a_2)V \\
S &@>>> &X_1\times X_2 &@>f_1\times f_2>>
&\M\times\M\end{CD}\]then
\[Z_S\isom S\times [\C^{d_1(S)}/\C^{d_0(S)}],\] where the vector space
$\C^{d_0(S)}$ acts trivially, and the map $u$ is the obvious
projection.
\end{prop}

\begin{pf}
By the existence of flattening stratifications, we can stratify
$X_1\times X_2$ by locally-closed subschemes $S$ such that the
restrictions of the families $E_1$ and $E_2$ have constant $\Ext$
groups. To ease the notation set $A=E_1$, $B=E_2$ and
\[V^k(S)=\lExt^k_{\OO_S}(B,A)\]
considered as a vector bundle over $S$.  The projection morphism
\[p\colon V^k(S)\to S\] defines an abelian group scheme over $S$. We will
show that \begin{equation} \label{object} Z_S\isom
[V^1(S)/V^0(S)],\end{equation} where the action of $V^0(S)$ on
$V^1(S)$ is the trivial one.  This will be enough because refining
the stratification if necessary we can assume that each of the
bundles $V^i(S)$ is trivial with fibre $\C^{d_i(S)}$. Then there
are obvious isomorphisms \[[V^1(S)/V^0(S)] \isom [S\times
\C^{d_1(S)}/ S \times \C^{d_0(S)}]\isom S \times
[\C^{d_1(S)}/\C^{d_0(S)}].\]

 Using Lemma \ref{headache} one
derives the following description of the stack $Z_S$. The objects
over a scheme $T$ consist of a morphism $f\colon T \to S$ and a
short exact sequence of $T$-flat sheaves on $T\times M$ of the
form
\begin{equation} \label{formm} 0\lra f^*(A) \lRa{\alpha} E
\lRa{\beta} f^*(B)\lra 0.\end{equation} Suppose given another such
map $g\colon U\to S$ and a sequence
 \[0\lra g^*(A) \lRa{\gamma} F \lRa{\delta} g^*(B)\lra 0.\]
Then a morphism between these two objects in $Z_S$ is a commuting
diagram of schemes
 \[ \xymatrix@C=.8em{
U\ar[dr]_{g}\ar[rr]^{h} && T\ar[dl]^{f}\\
&S }
\]
 and an isomorphism of
 sheaves $\theta\colon h^*(E) \to F$ such that the diagram
\begin{equation}\begin{CD}\label{c}
0 &@>>>& h^*f^*(A)& @>h^*(\alpha)>> &h^*(E) &@>h^*(\beta)>> &h^*f^*(B) &@>>> &0 \\
&&&&@V{\can}VV && @V{\theta}VV && @V{\can}VV \\
0& @>>>& g^*(A) &@>\gamma>> &F &@>\delta>> &g^*(B) &@>>>
&0\end{CD}\end{equation} commutes. Here $\can$ denotes the
canonical isomorphism. %\medskip

Lemma \ref{stanage} implies that there is a universal extension
class
\[\eta\in \Ext^1_{V^1(S)\times M} (p^*(B),p^*(A)).\]
Choose a
corresponding short exact sequence \begin{equation}
\label{morekids}0\lra p^*(A) \lRa{\gamma} F \lRa{\delta} p^*(B)
\lra 0.\end{equation} This defines an object of $Z_S(V^1(S))$ and
hence a morphism of stacks
\[q\colon V^1(S) \to Z_S.\]
%, such that the diagram
%\[ \xymatrix@C=.8em{
%V_1(T)\ar[dr]_{q}\ar[rr]^{q} && Z_T\ar[dl]^{f}\\
%&T }
%\]
%is 2-commutative.
%\medskip
Consider the fibre product \[W=V^1(S)\times_{Z_S}
V^1(S).\] For each scheme $T$ the groupoid $W(T)$ is a set, so we
can identify $W$ with the corresponding functor \[W \colon \sch/\C
\to \operatorname{Set}.\] It will be enough to show that the
functor $W$ is isomorphic to the functor
\[\Phi^0_{V^1(S)} (p^*(B),p^*(A))\]
and that there is a Cartesian diagram of stacks
\[\begin{CD} V^1(S)\times_S V^0(S) &@>\pi_1>>& V^1(S) \\
 @V\pi_1VV &&@VVqV \\
V^1(S) &@>q>> &Z_S\end{CD}\] Then $Z_S$ is  isomorphic to
the quotient stack corresponding to the trivial action of $V^0(S)$
on $V^1(S)$ as claimed.
%\medskip

 By the universal
property of $V^1(S)$ a morphism $a\colon T \to V^1(S)$ corresponds
to a morphism $f \colon T \to S$ together with an an extension
class
\[\zeta\in \Ext^1_{S\times M} (f^*(B),f^*(A)).\]
Under this correspondence $f=p\circ a$. The composite morphism
\[q\circ a\colon T \to Z_S\] then corresponds to the object of
$Z_S(T)$ defined by the morphism $f$ and the short exact sequence
\[0\lra f^*(A) \lra a^*(F) \lra f^*(B) \lra 0\]
obtained by applying $a^*$ to  \eqref{morekids} and composing with
the canonical isomorphisms. %\medskip

Suppose $b\colon T \to V^1(S)$ is another morphism corresponding
to a morphism $g\colon T \to S$ and an extension class \[\eta\in
\Ext^1_{S\times M} (g^*(B),g^*(A)).\] Then there is an isomorphism
of the corresponding objects of $Z_S(T)$ lying over the identity
of $T$ precisely if $f=g$ and there is an isomorphism of short
exact sequences
\begin{equation}\begin{CD}
0 &@>>>& f^*(A)& @>>> &a^*(F) &@>>> &f^*(B) &@>>> &0 \\
&&&&@V{=}VV && @VVV && @V{=}VV \\
0& @>>>& f^*(A) &@>>> &b^*(F) &@>>> &f^*(B) &@>>>
&0.\end{CD}\end{equation} In particular, it follows that
$\zeta=\eta$, and hence  by the universal property of $V^1(S)$ one
has $a=b$. Moreover the set of possible isomorphisms is in
bijection with
\[\Hom_{T\times M} (f^*(B),f^*(A) ).\]
Thus the elements of the set $W(T)$ consist of a morphism $a\colon
T \to V^1(S)$ and an element of
\[\Phi^0_{V^1(S)} (p^*(B),p^*(A))(a).\]
We leave it to the reader to check that this correspondence commutes with
pullback and hence defines an isomorphism of functors.
\end{pf}

% ******************************************************************************************************
% ******************************************************************************************************
% ******************************************************************************************************

\section{Proofs of Theorems \ref{fi} and \ref{se}}
Using Proposition \ref{tom} we can now give the proofs of Theorems \ref{fi} and
\ref{se}.

\subsection{Proof of Theorem \ref{fi}}
\label{firstproof}

Consider two elements \[a_i= [X_i\lRa{f_i} \M]\in \RHreg(\A),
\quad i=1,2,\] with $X_1$ and $X_2$ varieties. Let $E_i$ be the
family of coherent sheaves on $X_i$ corresponding to the map
$f_i$.  Stratify $X_1\times X_2$ by locally-closed subvarieties
$S_j$ as in Proposition \ref{tom}. In particular,  the vector
spaces
\[V^k(x_1,x_2)=\Ext^k_M\big(E_2|_{\{x_2\}\times M}, E_1|_{\{x_1\}\times M}\big),\qquad (x_1,x_2)\in S_j,\]
have constant dimension $d_k(S_j)$. Consider the diagram
\[\begin{CD}
Z_j &@>>> &Z &@>q>> \M^{(2)} @>h>> \M\\
@Vt_j VV && @VtVV  &@VV(a_1,a_2)V \\
S_j &@>>> &X_1\times X_2 &@>f_1\times f_2>> \M\times\M\end{CD}\]
According to Proposition \ref{tom} one has \[Z_j\isom
[Q_j/\C^{d_0(S_j)}]\] where $Q_j=V^1(S_j)$ is the total space of a
trivial vector bundle over $S_j$ with fibre $V^1(x_1,x_2)$ over a
point $(x_1,x_2)$. Since the $Z_j$ stratify $Z$ it follows that
\[a_1*a_2=[Z \lRa{b\circ h} \M]= \sum_j\, \LL^{-d_0(S_j)} [Q_j
\lRa{g_j} \M],\] which is regular. Here the morphism $g_j$ is
induced by the universal extension of the families $E_1$ and $E_2$
over $S_j$.
%\medskip

 For the second claim split $Q_j$ into the zero-section and its complement. The latter is a $\C^*$ bundle over the associated projective bundle,
 and it is easy to see that the morphism to $\M$ factors via this map. Thus \begin{equation}\label{domo} a_1*a_2=\sum_j
\LL^{-d_0(S_j)}\bigg([S_j\lRa{k}\M] + (\LL-1) [\PP( Q_j)\lRa{g_j}
\M]\bigg),\end{equation}where the morphism $k$ is induced by the
direct sum of the families  $E_1$ and $E_2$. We therefore obtain
\begin{equation}\label{bo}
a_1*a_2 = \sum_j [S_j\lRa{k} \M] = [X_1 \times X_2\lRa{k} \M] \mod
(\LL-1).\end{equation}   Clearly we would get the same answer if
we calculated $a_2*a_1$.\qed

% ******************************************************************************************************
% ******************************************************************************************************
% ******************************************************************************************************

\subsection{Proof of Theorem \ref{se}}
\label{secondproof}
We first check that the map $I$ is well-defined; it is then clearly unique. Stratify $\M$ by locally-closed substacks $\M_\tau$ such that $\lambda$ has constant value $\lambda(\tau)$ on $\M_\tau$.
There are projection maps \[\pi_\tau\colon \kst{\M}\to \kst{\M}\] defined by taking the fibre product with the inclusion $\M_\tau\subset \M$. For any $a\in\kst{\M}$
there is a canonical decomposition
\[a=\sum_{i} \pi_\tau(a),\]
where only finitely many of the terms are nonzero. If $a\in\kst{\M}$ is regular so are each of the $\pi_\tau(a)$.
 On the other hand if $b\in\kst{\M}$ is regular, we can project to an element of $\kst{\C}$,
  and using Lemma \ref{grimer} and Lemma \ref{toen} obtain a well-defined Euler characteristic $\chi(b)\in\Z$.
Thus we can define a group homomorphism $I$ by the formula
\[I(a)=\sum_i \lambda(\tau) \chi(\pi_\tau(a))\]
and this will clearly have the property stated in the Theorem.

Now take notation as in the proof of Theorem \ref{fi}. By Serre
duality, we have
\begin{equation}
\label{serre} V^k(x_1,x_2)=V^{3-k}(x_2,x_1)^*.\end{equation} Let
$\hat{Q}_j=V^2(S_j)$ be the bundle over $S_j$ whose fibre at
$(x_1,x_2)$ is $V^1(x_1,x_2)$. Let \[g_j\colon Q_j \to \M, \quad
\hat{g}_j\colon \hat{Q}_j\to \M,\]
 be the morphisms induced by taking the universal extensions of the families $E_1$ and $E_2$ over $S_j$.

We can assume that $f_i$ maps into $\M_{\alpha_i}\subset \M$ and
that $f_i^*(\lambda)$ is equal to the constant function with value
$n_i$. Then
\[I(a_i)=n_i \cdot \chi(X_i)\cdot x^{\alpha_i}.\]
Since $\chi(\LL)=1$, the expression \eqref{domo} shows that
\[I( a_1*a_2)=  \chi(X_1\times X_2,k^*(\lambda))\cdot x^{\alpha_1+\alpha_2}. \]
Using the first assumption \eqref{firstt}, we therefore obtain
\[ I(a_1* a_2)=\sigma^{\chi(\alpha_1,\alpha_2)}\cdot  n_1  n_2 \cdot \chi(X_1\times
X_2)= I(a_1)
* I(a_2).\]

To compute the Poisson bracket we  use \eqref{domo} again.
 Applying \eqref{serre} with $k=0$, and noting
that
\[\frac{\LL^{n}-\LL^{m}}{\LL-1}= n-m \mod (\LL-1),\]
we  obtain
\begin{align*}\{a_1,a_2\}=\sum_j \bigg(
(d_3(S_j)-d_0(S_j))&\cdot [S_j\lRa{k}\M]  \\&+  [\PP(
Q_j)\lRa{g_j}
\M]-[\PP(\hat{Q}_j)\lRa{\hat{g}_j}\M]\bigg).\end{align*} To
compute the Euler characteristic of a constructible function
 over $\PP(Q_j)$ it follows from \cite[Prop. 1]{mac} (see also \cite[Cor. 5.1]{verdier})
that we can first integrate over the fibres of the projection
\[\PP(Q_j)\to Q_j\] and then integrate the resulting constructible
function on the base $Q_j$. The second assumption \eqref{secondd}
together with $\chi(\PP(\C^n))=n$ therefore gives
 \[\chi(\PP(Q_j),g_j^*(\lambda)) - \chi(\PP(\hat{Q}_j),
 \hat{g}_j^*(\lambda))
 =d_1(S_j)\cdot\chi(S_j, k^*(\lambda)) -d_2(S_j)\cdot\chi(S_j,
 k^*(\lambda)),\]
and so
\begin{align*}I(\{a_1,a_2\})&=\chi(\alpha_1,\alpha_2)\cdot\chi(X_1\times X_2,k^*(\lambda))
 \cdot x^{\alpha_1+\alpha_2}\\
&=\sigma^{\chi(\alpha_1,\alpha_2)} \cdot n_1 n_2 \cdot
\chi(\alpha_1,\alpha_2)\cdot \chi(X_1\times X_2)\cdot
x^{\alpha_1+\alpha_2} =\{I(a_1),I(a_2)\}\end{align*}
 as required. \qed

% ********************************************************************************************
% ********************************************************************************************
% ********************************************************************************************

\begin{appendix}

\section{Fibre products of stacks}

%Given a stack $X$ and a scheme $S$  we denote by $X(S)$ the
%groupoid of morphisms $S\to X$.
Here we collect some well known material on fibre products of
stacks.

\subsection{Fibre product}
%The most important construction involving stacks for us will be the fibre product
Suppose given morphisms of stacks
\[f\colon X\to Z,\quad g\colon Y\to Z.\]
Recall the definition of the fibre product stack $X\times_Z Y$
 and the 2-commutative diagram
\[\begin{CD} X\times_Z Y &@>\pi_Y>> &Y \\
@V\pi_XVV && @VVgV \\
X &@>f>> &Z \\\end{CD}
\]
%\[\begin{CD} W &@>>> &Y \\
%@VVV && @VVgV \\
%X &@>>f> &Z \\\end{CD}%\medskip
%\]
The objects of $X\times_Z Y$  are triples $(x,y,\theta)$, where
$x$ and $y$ are objects of $X$ and $Y$ over the same scheme $S$,
and $\theta\colon f(x)\to g(y)$ is an isomorphism in the groupoid
$Z(S)$. A morphism
\[(\alpha,\beta)\colon (x,y,\theta)\to (x',y',\theta')\]
consists of morphisms $\alpha\colon x\to x'$ in $X$ and
$\beta\colon y\to y'$ in $Y$ such that the diagram
\[\begin{CD} f(x) &@>\theta>> &g(y) \\
@Vf(\alpha)VV && @VVg(\beta)V \\
f(x') &@>\theta'>> &g(y') \\\end{CD}%\medskip
\]
commutes. The morphisms $\pi_X$ and $\pi_Y$ are defined in the
obvious way.

We will call a morphism of stacks $f\colon X \to Z$ an
iso-fibration if the following property holds. Suppose $S$ is a
scheme and \[\theta\colon a \to b\] is an isomorphism in the
groupoid $Z(S)$. Suppose that there is an $a'\in X(S)$ such that
$f(a')=a$. Then there is an isomorphism \[\theta'\colon a'\to b'\]
in $X(S)$ such that $f(\theta')=\theta$. The following easy Lemma
will simplify many computations of such fibre products.

\begin{lemma}
\label{headache} With notation as above, define a full subcategory
$W\subset X\times_Z Y$ whose objects are triples $(x,y,\theta)$ as
above for which there is an object $z\in Z$ with
\[f(x)=z=g(y),  \quad \theta=\id_z.\]
 Suppose one of the morphisms $f$ or $g$ is an iso-fibration.
  Then  the inclusion functor $W\to X\times_Z Y$ is an equivalence of categories.
\end{lemma}

\begin{pf}
 For definiteness suppose that it is $f$ that is an iso-fibration. Given an object \[(x,y,\theta)\in X\times_Z Y\]
take an object $x'\in X$ such that $f(x')=g(y)$ and a morphism
$\alpha\colon x\to x'$ such that $f(\alpha)=\theta$. Then
\[(\alpha,\id_y)\colon (x,y,\theta)\to (x',y,\id)\] defines an
isomorphism. Thus $(x,y,\theta)$ is isomorphic to an object of the
subcategory $W$.
\end{pf}

\subsection{Cartesian diagrams} A diagram of stacks
\[\begin{CD} W &@>h>> &Y \\
@VjVV && @VVgV \\
X &@>f>> &Z \\\end{CD}
\]
is called 2-Cartesian if there is an equivalence of stacks
\[t\colon W \to X \times_Z Y\]
such that $j=\pi_X\circ t$ and $h=\pi_Y \circ t$. %An alternative
%definition would be that there is a 2-isomorphism \[\theta\colon
%g\circ h\isom f\circ j\] with a certain universal property.

\begin{lemma}
\label{nonse} Consider a 2-commutative diagram of the form
\[\begin{CD} V &@>>> &W &@>>> &Y\\
@VVV && @VVV && @VVV\\
U &@>>> &X &@>>> &Z\\\end{CD}
\]
and assume the right-hand small square to be 2-Cartesian. Then the
left hand small square is 2-Cartesian iff the big square is
2-Cartesian.
\end{lemma}

\begin{pf}
This is a standard fact, and we leave the proof to the reader.
\end{pf}

 We used the following
easy consequence many times in Section 4.

\begin{lemma}
\label{nonsense} Consider the following two diagrams of morphisms
of stacks
\[\begin{CD}
W  @>f>> {Y}   &\qquad\qquad\qquad&        {W}  @>f>> Y \\
@VgVV  @VVhV           @V(g,k\circ f)VV  @VV(h,k)V \\
X  @>j>> {Z}   &\qquad\qquad\qquad&     {{X} \times {T}}  @>(j,1)>> Z\times T \\
\end{CD}
\]
 Then if one is 2-Cartesian, so is the other.
\end{lemma}

\begin{pf}
This follows from Lemma \ref{nonse} once one knows that the square
\[\begin{CD} X\times T &@>(f,1)>> &Z\times T \\
@V\pi_XVV && @VV\pi_ZV \\
X &@>f>> &Z \\\end{CD}
\]
is 2-Cartesian. This follows from the diagram
\[\begin{CD} X \times T&@>f\times 1>> &Z\times T &@>\pi_T>> &T\\
@V\pi_XVV && @V\pi_ZVV && @VVV\\
X &@>f>> &Z&@>>> &\blob\\\end{CD}
\]by another application of Lemma
\ref{nonse}.
\end{pf}

\end{appendix}

\end{document}